\documentclass[12pt]{amsart}
\theoremstyle{plain}
\newtheorem{theorem}                 {Theorem}      [section]

\newtheorem{lemma}        [theorem]  {Lemma}

\theoremstyle{definition}
\newtheorem{example}      [theorem]  {Example}
\newtheorem{remark}       [theorem]  {Remark}

\DeclareMathOperator{\Rim}{Rim}

\numberwithin{equation}{section} \numberwithin{theorem}{section}

\title[]{\bf Natural Diagonal Riemannian Almost Product and Para-Hermitian Cotangent Bundles}
\author[]{S.~L.~Dru\c t\u a-Romaniuc\\ \\ \emph{The paper will appear in Czechoslovak Mathematical Journal}}

\begin{document}
\maketitle

\small \bf{Abstract:} \rm We obtain the natural diagonal almost
product and locally product structures on the total space of the
cotangent bundle of a Riemannian manifold. We find the Riemannian
almost product (locally product) and the (almost) para-Hermitian
cotangent bundles of natural diagonal lift type. We prove the
characterization theorem for the natural diagonal (almost)
para-K\"ahlerian structures on the total spaces of the cotangent
bundle.



\numberwithin{equation}{section}

\normalsize \section{Introduction} Some new interesting geometric
structures on the total space $T^*M$ of the cotangent bundle of a
Riemannian manifold $(M,g)$ were obtained for example in
\cite{Druta}, \cite{Munteanu}, \cite{OprPap1}--\cite{OprPor2} by
considering the natural lifts of the metric from the base manifold
to $T^*M$. Extensive literature, concerning the cotangent bundles
of natural bundles, may be found in \cite{Kolar1}.

The fundamental differences between the geometry of the cotangent
bundle and that of the tangent bundle, its dual, are due to the
construction of the lifts to $T^*M$, which is not similar to the
definition of the lifts to $TM$ (see \cite{YanoIsh}).

In a few papers such as \cite{An}-\cite{Cruceanu},
\cite{FZ}--\cite{Ivanov}, \cite{MihaiNic}, \cite{Peyghan}, and
\cite{Yano}, some almost product structures and almost
para-Hermitian structures (called also almost hyperbolic Hermitian
structures) were constructed on the total spaces of the tangent
and cotangent bundles.

In 1965, K. Yano initiated in \cite{Yano} the study of the
Riemannian almost product manifolds. A. M. Naveira gave in 1983 a
classification of these manifolds with respect to the covariant
derivative of the almost product structure (see \cite{Nav}). In
the paper \cite{SG} in 1992, M. Staikova and K. Gribachev obtained
a classification of the Riemannian almost product manifolds, for
which the trace of the almost product structure vanishes, the
basic class being that of the almost product manifolds with
nonintegrable structure (see \cite{Mek}).

A classification of the almost para-Hermitian manifolds was made
in 1988 by C. Bejan, who obtained in \cite{Be} 36 classes, up to
duality, and the characterizations of some of them. P. M. Gadea
and J. Mu\~noz Masqu\'e gave in 1991 a classification \`{a} la
Gray-Hervella, obtaining 136 classes, up to duality (see
\cite{GM}). Maybe the best known class of (almost) para-Hermitian
manifolds are the (almost) para-K\"ahler manifolds, characterized
by the closure of the associated 2-form, and studied for example
in \cite{AMT} and \cite{Dorota}.

In the present paper we consider an $(1,1)$-tensor field $P$
obtained as a natural diagonal lift of the metric $g$ from the
base manifold $M$ to the total space $T^*M$ of the cotangent
bundle. This tensor field depends on four coefficients which are
smooth functions of the energy density $t$. We first determine the
conditions under which the tensor field constructed in this way is
an almost product structure on $T^*M$. We obtain some simple
relations between the coefficients of $P$. From the study of the
integrability conditions of the determined almost product
structure, it follows that the base manifold must be a space form,
and two coefficients may be expressed as simple rational functions
of the other two coefficients, their first order derivatives, the
energy density, and the constant sectional curvature of the base
manifold. Then we prove the characterization theorems for the
cotangent bundles which are Riemannian almost product (locally
product) manifolds, or (almost) para-Hermitian manifolds, with
respect to the obtained almost product structure, and a natural
diagonal lifted metric $G$. Finally we obtain the (almost)
para-K\"ahler cotangent bundles of natural diagonal lift type.

Throughout this paper, the manifolds, tensor fields and other
geometric objects are assumed to be differentiable of class
$C^\infty$ (i.e. smooth). The Einstein summation convention is
used , the range of the indices $h,i,j,k,l,m,r, $ being always
$\{1,\dots ,n\}$.

\section{Preliminary results}
The cotangent bundle of a smooth $n$-dimensional Riemannian
manifold may be endowed with a structure of $2n$-dimensional
smooth manifold, induced by the structure on the base manifold. If
$(M,g)$ is a smooth Riemannian manifold of dimension $n$, we
denote its cotangent bundle by $\pi :T^*M\rightarrow M$. Every
local chart on $M$, $(U,\varphi )=(U,x^1,\dots ,x^n)$ induces a
local chart $(\pi^{-1}(U),\Phi )=(\pi^{-1}(U),q^1,\dots , q^n,$
$p_1,\dots ,p_n)$, on $T^*M$, as follows. For a cotangent vector
$p\in \pi^{-1}(U)\subset T^*M$, the first $n$ local coordinates
$q^1,\dots ,q^n$ are  the local coordinates of its base point
$x=\pi (p)$ in the local chart $(U,\varphi )$ (in fact we have
$q^i=\pi ^* x^i=x^i\circ \pi, \ i=1,\dots n)$. The last $n$ local
coordinates $p_1,\dots ,p_n$ of $p\in \pi^{-1}(U)$ are the vector
space coordinates of $p$ with respect to the basis
$(dx^1_{\pi(p)},\dots , dx^n_{\pi(p)})$, defined by the local
chart $(U,\varphi )$,\ i.e. $p=p_idx^i_{\pi(p)}$.

The concept of $M$-tensor field on the cotangent bundle of a
Riemannian manifold was defined by the present author in
\cite{Druta}, in the same manner as the $M$-tensor fields were
introduced on the tangent bundle (see \cite{Mok}).

We recall the splitting of the tangent bundle to $T^*M$ into the
vertical distribution $VT^*M= {\rm Ker}\ \pi _*$ and the
horizontal one determined by the Levi Civita connection $\dot
\nabla $ of the metric $g$:
\begin{eqnarray}\label{descomp}
~~~~~~~~~~~~~~~~~~~~~~~~~~TT^*M=VT^*M\oplus HT^*M.
\end{eqnarray}
If $(\pi^{-1}(U),\Phi)=(\pi^{-1}(U),q^1,\dots ,q^n,p_1,\dots
,p_n)$ is a local chart on $T^*M$, induced from the local chart
$(U,\varphi )= (U,x^1,\dots ,x^n)$, the local vector fields
$\frac{\partial}{\partial p_1}, \dots , \frac{\partial}{\partial
p_n}$ on $\pi^{-1}(U)$ define a local frame for $VT^*M$ over $\pi
^{-1}(U)$ and the local vector fields $\frac{\delta}{\delta
q^1},\dots ,\frac{\delta}{\delta q^n}$ define a local frame for
$HT^*M$ over $\pi^{-1}(U)$, where $\frac{\delta}{\delta
q^i}=\frac{\partial}{\partial q^i}+\Gamma^0_{ih}
\frac{\partial}{\partial p_h},\ \Gamma ^0_{ih}=p_k\Gamma ^k_{ih},$
and $\Gamma ^k_{ih}(\pi(p))$ are the Christoffel symbols of $g$.

The set of vector fields $\{\frac{\partial}{\partial
p_i},\frac{\delta}{\delta q^j}\}_{i,j=\overline{1,n}}$, denoted by
$\{\partial^i, \delta_j\}_{i,j=\overline{1,n}}$, defines a local
frame on $T^*M$, adapted to the direct sum decomposition
(\ref{descomp}).

We consider the energy density defined by $g$ in the cotangent
vector $p$:
\begin{eqnarray*}
\begin{array}{c}
t=\frac{1}{2}\|p\|^2=\frac{1}{2}g^{-1}_{\pi(p)}(p,p)=\frac{1}{2}g^{ik}(x)p_ip_k,
\ \ \ p\in \pi^{-1}(U).
\end{array}
\end{eqnarray*}

We have $t\in [0,\infty)$ for all $p\in T^*M$.

In the sequel we shall use the following lemma, which may be
proved easily.

\begin{lemma}\label{lema1}
\it{If $n>1$ and $u,v$ are smooth functions on $T^*M$ such that
\begin{eqnarray*}
u g_{ij}+v p_ip_j=0, \quad u g^{ij}+v g^{0i}g^{0j}=0,\quad or \
u\delta ^i_j+vg^{0i} p_j=0,
\end{eqnarray*}
on the domain of any induced local chart on $T^*M$, then $u=0,\
v=0$.} We used the notation $g^{0i}=p_hg^{hi}$.
\end{lemma}

\section{Almost product structures of natural diagonal lift type on the cotangent bundle}
In this section we shall find the almost product structures on the
(total space of the) cotangent bundle, which are natural diagonal
lifts of the metric from the base manifold $M$ to $T^*M$. Then we
shall study the integrability conditions for the determined
structures, obtaining the natural diagonal locally product
structures on $T^*M$.

An \emph{almost product structure} $J$ on a differentiable
manifold $M$ is an $(1,1)$- tensor field on $M$ such that $J^2=I$.
The pair $(M,J)$ is called an \emph{almost product manifold}. When
the almost product structure $J$ is integrable, it is called a
\emph{locally product structure}, and the manifold $(M,J)$ is a
\emph{locally product manifold}.

An \emph{almost paracomplex manifold} is an almost product
manifold $(M,J)$, such that the two eigenbundles associated to the
two eigenvalues $+1$ and $-1$ of $J$, respectively, have the same
rank. Equivalently, a splitting of the tangent bundle $TM$ into
the Whitney sum of two subbundles $T^\pm M$ of the same fiber
dimension is called an \emph{almost paracomplex structure} on $M$.

V. Cruceanu presented in \cite{Cruceanu} two simple almost product
structures on the total space $T^*M$ of the cotangent bundle,
obtained by considering on the base manifold $M$ a linear
connection $\nabla$ and a non-degenerate $(0,2)$-tensor field $g$.
If $\alpha$ is a differentiable 1-form and $X$ is a vector field
on $M$, $\alpha^V$ denotes the vertical lift of $\alpha$ and $X^H$
the horizontal lift of $X$ to $T^*M$, one can consider
\begin{equation}\label{P}
~~~~~~~~~~~~~~~~~~~~~~~~~~P(X^H)=-X^H,\ P(\alpha^V)=\alpha^V,
\end{equation}
\begin{equation}\label{Q}
~~~~~~~~~~~~~~~~~~~~~~~~~~Q(X^H)=(X^\flat)^V,\
Q(\alpha^V)=(\alpha^\sharp)^H,
\end{equation}
where $X^\flat=g_X$ is the 1-form on $M$ defined by
$X^\flat(Y)=g_X(Y)=g(X,Y), \forall Y\in \mathcal{T}^1_0(M)$,
$\alpha^\sharp=g^{-1}_\alpha$ is a vector field on $M$ defined by
$g(\alpha^\sharp,Y)=\alpha (Y),~\forall Y \in \mathcal{T}^1_0(M)$.
$P$ is a paracomplex structure if and only if $\nabla$ has
vanishing curvature, while $Q$ is paracomplex if and only if the
curvature of $\nabla$ and the exterior covariant differential $Dg$
of $g$, given by
$$
(Dg)(X,Y)=\nabla_X(Y^\flat)-\nabla_Y(X^\flat)-[X,Y]^\flat,
$$
vanishes.

The results from \cite{Kolar} and \cite{KowalskiSek} concerning
the natural lifts, allow us to introduce an $(1,1)$ tensor field
$P$ on $T^*M$, which is a natural diagonal lift of the metric $g$
from the base manifold to the total space $T^*M$ of cotangent
bundle. Using the adapted frame
$\{\partial^i,\delta_j,\}_{i,j=\overline{1,n}}$ to $T^*M$, we
define $P$ by the relations:
\begin{eqnarray}\label{defPcomp}
P\delta_i=P^{(1)}_{ij}\partial^j, \quad
P\partial^i=P_{(2)}^{ij}\delta_j,
\end{eqnarray}
where the $M-$tensor fields involved as coefficients have the
forms
\begin{eqnarray}\label{compP}
~~~~~~~P^{(1)}_{ij}=a_1(t)g_{ij}+b_1(t)p_ip_j,\quad
P_{(2)}^{ij}=a_2(t)g^{ij}+b_2(t)g^{0i}g^{0j}.
\end{eqnarray}
$a_1,\ b_1,\ a_2,$ and $b_2$ being smooth functions of the energy
density $t$.

The  invariant expression of the defined structure is
\begin{eqnarray}\label{defP}
\begin{array}{c}PX^H_p=a_1(t)
(X^\flat)^V_p+b_1(t)p(X)p_p^V, \\\\
P\alpha^V_p=a_2(t)(\alpha^\sharp)_p^H+b_2(t)g^{-1}_{\pi(p)}
(p,\alpha)(p^\sharp)_p^H,
\end{array}
\end{eqnarray}
in every point $p$ of the induced local card $(\pi^{-1}(U),\Phi)$
on $T^*M$, $\forall ~X \in \mathcal{T}^1_0(M), \forall~ \alpha \in
\mathcal{T}^0_1(M)$. The vector $p^\sharp$ is tangent to $M$ in
$\pi (p)$, $p^V=p_i\partial^i$ is the Liouville vector field on
$T^*M$ , and $(p^\sharp)^H=g^{0i}\delta_i$ is the geodesic spray
on $T^*M$.

\normalsize

\begin{example}\label{ex1}
When $a_1=a_2=1,$ $b_1$ and $b_2$ vanish, we have the structure
given by $(\ref{Q}).$
\end{example}

The following theorems present the conditions under which the
above tensor field $P$ defines an almost product (locally product)
structure on the total space of the cotangent bundle.

\begin{theorem}\label{th4}
{The tensor field $P$, given by \rm (\ref{defPcomp}) \it or
\rm(\ref{defP}), \it defines an almost product structure of
natural diagonal lift type on $T^*M$, if and only if its
coefficients satisfy the relations
\begin{equation}\label{inlocuire}
~~~~~~~~~~~~~~~~~~~~a_1=\frac{1}{a_2}\ ,\ \ \
a_1+2tb_1=\frac{1}{a_2+2tb_2}.
\end{equation}}
\end{theorem}

\proof
The condition $P^2=I$ in the definition of the almost
product structure, may be written in the following form, by using
(\ref{defPcomp}):
$$
P^{(1)}_{ij}P_{(2)}^{il}=\delta^l_j,~~~~~~
P_{(2)}^{ij}P^{(1)}_{il}=\delta^j_l,
$$
and replacing (\ref{compP}) it becomes
$$
(a_1 a_2-1)\delta^l_j + [b_1 (a_2 + 2 tb_2)+ a_1 b_2] g^{0l}p_j=0.
$$

Using Lemma \ref{lema1}, we have that the above expression vanish
if and only if
$$
a_1=\frac{1}{a_2},\ b_1=-\frac{a_1 b_2}{a_2 + 2 tb_2},
$$
which imply also the second relation in (\ref{inlocuire}).
\endproof

\begin{remark}\label{remark1}
If $a_1=\frac{1}{\beta},\ a_2=\beta$, $b_1=\frac{u}{\alpha\beta}$
and $b_2=-\frac{u\beta}{\alpha +2 t u}$, where $\alpha$ and
$\beta$ are real constants and $u$ is a smooth function of $t$,
the statements of \rm Theorem \ref{th4} \it are satisfied, so the
structure considered in \rm \cite{Peyghan} \it is an almost
product structure on the total space $T^*M$ of the cotangent
bundle.
\end{remark}

\begin{theorem}\label{th3}
The natural diagonal almost product structure $P$ on the total
space of the cotangent bundle of an $n$-dimensional connected
Riemannian manifold $(M,g)$, with $n>2$, is a locally product
structure on $T^*M$ (i.e. $P$ is integrable) if and only if the
base manifold is of constant sectional curvature $c$, and the
coefficients $b_1,\ b_2$ are given by:
\begin{equation}\label{b}
~~~~~~~~~~~~~~~~~~~~~~b_1=\frac{a_1 a_1^\prime + c}{a_1 - 2t
a_1^\prime},\ b_2=\frac{a_1 a_2^\prime - a_2^2 c}{a_1 + 2 c t a_2
}.
\end{equation}
\end{theorem}

\proof
The almost product structure $P$ on $T^*M$ is integrable if
and only if the vanishing condition for the Nijenhuis tensor field
$N_P$,
$$
N_P(X,Y)=[PX,PY]-P[PX,Y]-P[X,PY]+P^2[X,Y],\forall X,Y\in
\mathcal{T}^1_0(T^*M),
$$
is satisfied.

Studying the components of $N_P$ with respect to the adapted frame
on $T^*M$, $\{\partial^i,\delta_j\}_{i,j=\overline{1,n}}$, we
first obtain:
$$
N_P(\partial^i,\partial^j)=[P^{(1)}_{km}(\partial^jP_{(2)}^{mi}-\partial^iP_{(2)}^{mj})+
\Rim^0_{kml} P_{(2)}^{mi}P_{(2)}^{lj}]\partial^k.
$$

Replacing the values (\ref{compP}), after some tensorial
computations the above expression becomes
\begin{equation}\label{curvature}
a_1(a_2^\prime - b_2)(\delta^h_ip_j-\delta^h_jp_i)
-a_2^2\Rim^h_{kij} +
a_2b_2(\Rim^h_{kjl}p_i-\Rim^h_{kil}p_j)g^{0k}g^{0l}=0.
\end{equation}

We differentiate (\ref{curvature}) with respect to $p_h$. Since
the curvature of the base manifold does not depend on $p$, we take
the value of this derivative  at $p=0$, and we obtain
\begin{equation}\label{sectional curvature}
~~~~~~~~~~~~~~~~~~~~~~~~~~~R^h_{kij}=c(\delta^h_ig_{kj}-\delta^h_jg_{ki}),
\end{equation}
where
$$
c=\frac{a_1(0)}{a_2^2(0)}(a'_2(0)-b_2(0))
$$
is a function depending on $q^1,...,q^n$ only. Schur's theorem
implies that $c$ must be a constant when $M$ is connected, of
dimension $n>2$.

Moreover, by using (\ref{sectional curvature}), the relation
(\ref{curvature}) becomes
\begin{equation}\label{rel}
~~~~~~~~~~[a_1 a_2^\prime - a_2^2 c - b_2 (a_1 + 2 t a_2 c)]
   (\delta^h_ig_{kj}-\delta^h_jg_{ki})=0.
\end{equation}

Solving (\ref{rel}) with respect to $b_2$, we obtain the second
relation in (\ref{b}).

The Nijenhuis tensor field computed for both horizontal arguments
is
$$
N_P(\delta_i,\delta_j)=(P^{(1)}_{li} \partial^lP^{(1)}_{hj} -
P^{(1)}_{lj}
\partial^lP^{(1)}_{hi}+ \Rim^0_{hij})\partial^h,
$$
which vanishes if and only if
$$
[b_1(2 t a_1^\prime  - a_1)+a_1^\prime a_1 + c]
(g_{hj}p_i-g_{hi}p_j)=0,
$$
namely when $b_1$ has the expression in (\ref{b}).

The mixed components of the Nijenhuis tensor field, have the forms
$$
N_P(\delta_i,\partial^j)=-N_P(\partial^j,\delta_i)=(P^{(1)}_{mi}
\partial^m P_{(2)}^{hj} + P_{(2)}^{hl} \partial^j P^{(1)}_{li} -
P_{(2)}^{lh}P_{(2)}^{jm} \Rim^0_{lim})\delta_h,
$$
which after replacing (\ref{compP}) and (\ref{sectional
curvature}) become
$$
(a_1 a_2' + a_2 b_1 - a_2^2 c + 2t a_2' b_1) g^{hj} p_i
$$
$$
+ (a_2 b_1 + a_1 b_2 + 2 t b_1 b_2) \delta^j_i g^{0h} + (a_1' a_2
+ a_1 b_2 + a_2^2 c + 2 cta_2 b_2) \delta^h_i g^{0j}
$$
$$+
 (a_2 b_1' + a_1' b_2 + 3 b_1 b_2 + a_1 b_2' - a_2 b_2 c +
    2t b_1' b_2 + 2t b_1 b_2') p_i g^{0h}g^{0j}.
$$

Taking (\ref{inlocuire}) into account, the above expression takes
the form
$$
\frac{(a_1- 2a_1't)b_1-a_1a_1'- c }{a_1^2}g^{hj}p_i +
\frac{a_1a_1'+ c - (a_1- 2a_1't)b_1}{a_1(a_1 + 2tb_1)}
    \delta^h_ig^{0j}
$$
$$
+ \frac{(a_1a_1'+ c)b_1 - (a_1 - 2ta_1')b_1^2  }{a_1^2(a_1 +
2tb_1)}p_ig^{0h}g^{0j}
$$
and it vanishes if and only if $b_1$ is expressed by the first
relation in (\ref{b}).

One can verify that all the components of the Nijenhuis tensor
field vanish under the same conditions, so the almost product
structure $P$ on $T^*M$ is integrable, i.e. $P$ is a locally
product structure on $T^*M$.
\endproof

\begin{remark}\label{remark3}
If the coefficients involved in the definition of $P$ have the
values presented in \rm Remark \ref{remark1}, \it the relations
\rm (\ref{b}) \it take the form $u=c\alpha\beta^2$, so \rm Theorem
\ref{th3} \it implies the results stated in \rm{\cite[Theorem
4.2]{Peyghan}}.
\end{remark}

\section{Natural diagonal Riemannian almost product and almost para-Hermitian
structures on $T^*M$} Authors like M. Anastasiei, C. Bejan, V.
Cruceanu, H. Farran, A. Heydari, S. Ishihara, I. Mihai, G. Mitric,
C. Nicolau, V. Oproiu, L. Ornea, N. Papaghiuc, E. Peyghan, K.
Yano, and M.~S. Zanoun considered almost product structures and
almost para-Hermitian structures (called also almost hyperbolic
Hermitian structures) on the total spaces of the tangent and
cotangent bundles.

A Riemannian manifold $(M,g)$, endowed with an almost product
structure $J$, satisfying the relation
\begin{equation}\label{defprodpara}
~~~~~~~~~~~~~~~~~~~g(JX,JY)=\varepsilon g(X,Y),\ \forall X,Y\in
\mathcal{T}^1_0(M),
\end{equation}
is called a \emph{Riemannian almost product manifold} if
$\varepsilon=1$, or an \emph{almost para-Hermitian manifold}
(called also an \emph{almost hyperbolic Hermitian manifold}) if
$\varepsilon=-1$.

In the following we shall find the Riemannian almost product
(locally product) and the (almost) para-Hermitian cotangent
bundles of natural diagonal lift type.

To this aim, we consider a natural diagonal lifted metric on the
total space $T^*M$ of the cotangent bundle, defined by:
\begin{equation}\label{Ginvar}
\left\{
\begin{array}{l}
G_p(X^H, Y^H) = c_1(t)g_{\pi(p)}(X,Y) + d_1(t)p(X)p(Y),
\\\\
G_p(\alpha^V,\omega^V) = c_2(t)g^{-1}_{\pi(p)}(\alpha,\omega) +
d_2(t)g^{-1}_{\pi(p)}(p,\alpha)g^{-1}_{\pi(p)}(p,\omega),
\\\\
G_p(X^H,\alpha^V) = G_p(\alpha^V,X^H) =0,
\end{array}
\right.
\end{equation}
$\forall~ X,Y \in \mathcal{T}^1_0(M),$ $\forall~ \alpha, \omega
\in \mathcal{T}^0_1(M), \forall~p \in T^*M$, where the
coefficients $c_1,\ c_2,\ d_1,\ d_2$ are smooth functions of the
energy density.

The conditions for $G$ to be nondegenerate are assured if
$$
c_1c_2\neq0, \ (c_1+2td_1)(c_2+2td_2)\neq 0.
$$

The metric $G$ is positive definite if
$$
c_1+2td_1>0,\quad c_2+2td_2>0.
$$

Using the adapted frame $\{\partial^i,
\delta_j\}_{i,j=\overline{1,n}}$ on $T^*M$, (\ref{Ginvar}) becomes
\begin{eqnarray}\label{compG}
\begin{cases}
G(\delta_i, \delta_j)=G^{(1)}_{ij}=c_1(t)g_{ij}+ d_1(t)p_ip_j,\\
G(\partial^i,\partial^j)=G_{(2)}^{ij}=c_2(t)g^{ij}+d_2(t)g^{0i}g^{0j},
\\
G(\partial^i,\delta_j)=G(\delta_i,\partial^j)=0,
\end{cases}
\end{eqnarray}
where $c_1,\ c_2,\ d_1,\ d_2$ are smooth functions of the density
energy on $T^*M$.

Next we shall prove the following characterization theorem:

\begin{theorem}\label{aprprod}
Let $(M,g)$ be an $n$-dimensional connected Riemannian manifold,
with $n>2$, and $T^*M$ the total space of its cotangent bundle.
Let $G$ be a natural diagonal lifted metric on $T^*M$, defined by
\rm(\ref{Ginvar}), \it and $P$ an almost product structure on
$T^*M$, characterized by \rm Theorem $\ref{th4}$. \it Then
$(T^*M,G,P)$ is a Riemannian almost product manifold, or an almost
para-Hermitian manifold if and only if the following
proportionality relations between the coefficients hold:
\begin{eqnarray}\label{almprod1}
~~~~~~~~~~~~\frac{c_1}{a_1}=\varepsilon\frac{c_2}{a_2}=\lambda, \
\frac{c_1+2td_1}{a_1+2tb_1}=\varepsilon\frac{c_2+2td_2}{a_2+2tb_2}=\lambda
+2 t \mu,
\end{eqnarray}
where $\varepsilon$ takes the corresponding values from the
definition $(\ref{defprodpara})$, and the proportionality
coefficients $\lambda>0 $ and $\lambda +2t\mu>0$ are some
functions depending on the energy density $t$.

If moreover, the relations stated in \rm Theorem \ref{th3} \it are
fulfilled, then $(T^*M,G,$ $P)$ is a Riemannian locally product
manifold for $\varepsilon=1$, or a para-Hermitian manifold for
$\varepsilon=-1$.
\end{theorem}

\proof With respect to the adapted frame
$\{\partial^i,\delta_j\}_{i,j=\overline{1,n}}$, the relation
(\ref{defprodpara}) becomes:
\begin{equation}\label{condalmprod}
~~~~G(P\delta_i,P\delta_j)=\varepsilon G(\delta_i,\delta_j),\
G(P\partial^i,P\partial^j)=\varepsilon G(\partial^i,\partial^j),\
G(P\partial^i,P\delta_j)=0,
\end{equation}
and using (\ref{defPcomp}) and (\ref{compG}) we have
\begin{eqnarray*}
~~~~~~~~~~~~~~\begin{array}{c} (-\varepsilon c_1 + a_1^2 c_2)
g_{ij} + [ -\varepsilon d_1 + a_1^2 d_2+2  b_1 c_2(a_1 +t
 b_1) +
    4 t b_1 d_2  (a_1+ tb_1 )]p_ip_j=0,
\\\\
(a_2^2 c_1 - \varepsilon c_2) g^{ij} + [-\varepsilon d_2 + a_2^2
d_1 + 2 b_2 c_1 (a_2 + t b_2) +
    4 t b_2 d_1 (a_2 + tb_2 )] g^{0i}g^{0j}=0.
\end{array}
\end{eqnarray*}

Taking into account Lemma \ref{lema1}, the coefficients which
appear in the above expressions vanish. Due to the first relation
in (\ref{inlocuire}), we get by equalizing to zero the
coefficients of $g_{ij}$ and $g^{ij}$, the first relation in
(\ref{almprod1}).

Moreover, multiplying by $2t$ the coefficients of $p_ip_j$ and
$g^{0i}g^{0j}$ and adding them to the coefficients of $g_{ij}$ and
$g^{ij}$, respectively, we obtain:
\begin{eqnarray}\label{sist}
~~~~~~~~~~~~~~\begin{array}{l} -\varepsilon(c_1+2td_1) +
(a_1+2tb_1)^2 (c_2+2td_2)=0,
\\\\
(a_2+2tb_2)^2 (c_1+2tb_1) -\varepsilon (c_2+2td_2)=0.
\end{array}
\end{eqnarray}
Using the second relation in (\ref{inlocuire}), (\ref{sist}) leads
to the second relation in (\ref{almprod1}).
\endproof

\begin{remark} \label{remark2}When the coefficients of the almost product
structure $P$ have the expressions in \rm Remark \ref{remark1},
\it and the coefficients of the metric $G$ on $T^*M$ are
$c_1=a_1,$ $d_1=b_1,$ $c_2=-a_2$ and $d_2=-b_2,$ \rm Theorem
\ref{aprprod} \it implies that $T^*M$ endowed with the almost
product structure and with the metric considered in \rm
\cite{Peyghan} \it is an almost para-Hermitian manifold.
\end{remark}

\section{Natural diagonal para-K\"ahler structures on $T^*M$}

In the sequel we shall study the cotangent bundles endowed with
para-K\"ahler structures of natural diagonal lift type. This class
of almost para-Hermitian structures, studied for example in
\cite{AMT} and \cite{Dorota}, is characterized by the closure of
the associated 2-form $\Omega $.

The 2-form $\Omega $ associated to the almost para-Hermitian
structure $(G,P)$ of natural diagonal lift type on the total space
of the cotangent bundle is given by the relation
$$
\Omega (X,Y)=G(X,PY), \forall X,Y\in \mathcal{T}^1_0(T^*M).
$$

Studying the closure of $\Omega$, we may prove the following
theorem:

\begin{theorem}\label{th6}
The almost para-Hermitian structure $(G,P)$ of natural diagonal
lift type on the total space $T^*M$ of the cotangent bundle of a
Riemannian manifold $(M,g)$ is almost para-K\"{a}hlerian if and
only if
$$
\mu=\lambda ^\prime .
$$
\end{theorem}

\proof \rm The 2-form $\Omega$ on $T^*M$ has the following
expression with respect to the local adapted frame
$\{\partial^i,\delta_j\}_{i,j=1,\dots,n}$:
$$
\Omega(\partial^i,\partial^j)=\Omega(\delta_i,\delta_j)=0,\quad
\Omega(\partial^j,\delta_i)=G_{(2)}^{jh}P^{(1)}_{hi},\quad
\Omega(\delta_i,\partial^j)=G^{(1)}_{ih}P_{(2)}^{hj},
$$

By substituting (\ref{compP}) and (\ref{compG}) in the above
expressions, and taking into account the conditions for
$(T^*M,G,P)$ to be an almost para-Hermitian manifold (see Theorem
\ref{aprprod}), we have
\begin{equation}\label{omega}
~~~~~~~~~~~~~~~~~~~\Omega(\delta_i,\partial^j)=-\Omega(\partial^j,\delta_i)
=\lambda \delta^j_i+\mu p_ig^{0j},
\end{equation}
which has the invariant expression
$$
\Omega\left(X_p^H,\alpha_p^V\right)= \lambda \alpha(X)+\mu
p(X)g^{-1}_{\pi(p)}(p,\alpha),
$$
for every  $X\in \mathcal{T}^1_0(M),\ \alpha\in
\mathcal{T}^0_1(M),\ p \in T^*M,$

Taking into account the relation (\ref{omega}) we have that the
$2$-form $\Omega$ associated to the natural diagonal
para-Hermitian structure has the form
\begin{equation}
~~~~~~~~~~~~~~~~~~~~~~~~\Omega =(\lambda \delta^j_i+\mu p_i
g^{0j})dq^i\wedge Dp_j,
\end{equation}
where  $Dp_j=dp_j-\Gamma^0_{jh}dq^h$ is the absolute differential
of $p_i$.

Moreover, the differential of $\Omega$ will be
$$
d\Omega=(d\lambda \delta^j_i+d\mu g^{0j}p_i+\mu d g^{0j}p_i+\mu
g^{0j}dp_i)\wedge dq^i\wedge Dp_j  - (\lambda \delta^j_i+\mu
p_ig^{0j})dq^i\wedge dDp_j .
$$

Let us compute the expressions of $d\lambda,\ d\mu, \ dg^{0i}$ and
$dDp_i$:
$$
d\lambda=\lambda ' g^{0h}Dp_h,\quad d\mu=\mu ' g^{0h}Dp_h,\quad
dg^{0i}=g^{hi}Dp_h-\Gamma^i_{j0}dq^j,
$$
$$
dDp_i=\frac{1}{2}R^0_{ijh}dq^h\wedge dq^j-\Gamma^h_{ij}Dp_h\wedge
dq^j.
$$

Then, by substituting these relations into the expression of
$d\Omega$, taking into account the properties of the external
product, the symmetry of $g^{ij}$ and $\Gamma^h_{ij}$ and the
Bianchi identities, we obtain
$$
d\Omega=(\mu-\lambda^\prime)p_kg^{kh}\delta^i_jDp_h\wedge
Dp_i\wedge dq^j,
$$
which, due to the antisymmetry of $\delta^i_jDp_i\wedge dq^j$, may
be written as
$$
d\Omega=\frac{1}{2}(\mu-\lambda^\prime)p_k(g^{kh}\delta^j_i-g^{kj}\delta^h_i)Dp_h\wedge
Dp_j\wedge dq^i,
$$
and it vanishes if and only if $\mu=\lambda '$. \endproof

\vskip 2mm

Using the theorems \ref{th4}, \ref{th3} and \ref{th6}, we
immediately prove:

\begin{theorem}\label{paraK}
An almost para-Hermitian structure $(G,P)$ of natural diagonal
lift type on the total space $T^*M$ of the cotangent bundle of a
Riemannian manifold $(M,g)$ is para-K\"ahlerian if and only if $P$
is a locally product structure (see \rm Theorem $\ref{th3}$\it)
and $\mu=\lambda^\prime$.
\end{theorem}

\begin{remark}
The almost para-K\"ahlerian structures of natural diagonal lift
type on $T^*M$ depend on three essential coefficients $a_1, \ b_1$
and $\lambda$, while the natural diagonal para-K\"ahlerian
structures on $T^*M$ depend on two essential coefficients $a_1$
and $\lambda$, which in both cases must satisfy the supplementary
conditions $a_1>0,\ a_1+2tb_1>0,\ \lambda>0, \ \lambda +2 t
\lambda'>0$, where $b_1$ is given by $(\ref{b}).$
\end{remark}

\begin{remark}
Taking into account \rm Remark \ref{remark2}, \it we have that the
constant $\lambda$ is equal to $1$ and $\mu$ vanishes, so \rm
Theorem \ref{th6} \it leads to the statements in \rm{\cite[Theorem
3.1]{Peyghan}}, \it namely the structure constructed in \rm
\cite{Peyghan} \it is almost para-K\"ahlerian on $T^*M$. Moreover,
taking into account \rm Remark \ref{remark3}, \it it follows that
the relations \rm (14) \it in \rm \cite{Peyghan} \it are fulfilled
in the case when the constructed structure is para-K\"ahler.
\end{remark}

\textbf{Acknowledgements.} The author would like to thank
Professor V. Oproiu for the techniques learned and Professor M. I.
Munteanu for the critical remarks during the preparation of this
work.

The paper was supported by the Program POSDRU/89/1.5/S/49944, "Al.
I. Cuza" University of Ia\c si, Romania.

{\small {\em Author's address}:
 {\em Simona-Luiza Dru\c t\u
a-Romaniuc}, Department of Sciences, 54 Lascar Catargi Street
RO-700107,  "Al. I. Cuza" University of Ia\c si, Romania, e-mail:
simonadruta@yahoo.com.}

\end{document}